\documentclass[a4paper,twoside,12pt]{article}
\usepackage[english]{babel}
\usepackage{graphicx}
\usepackage{bbding}
\usepackage[T1]{fontenc}
\usepackage{kpfonts}
\usepackage{fancyhdr}
\usepackage{amsmath,amssymb,color,epsfig}
\usepackage{mathrsfs}
\usepackage{eufrak}
\usepackage{dsfont}
\usepackage{upgreek}
\usepackage{fancyhdr}

\newtheorem{theorem}{Theorem}[section]
\newtheorem{proposition}[theorem]{Proposition}
\newtheorem{corollary}[theorem]{Corollary}

\newtheorem{lemma}[theorem]{Lemma}
\newtheorem{definition}[theorem]{Definition}

\def\tr{{\rm tr}}

\def\1{\mathds{1}}
\def\e{\varepsilon}

\title{Generalized Gaussian Random Unitary Matrices Ensemble}
\date{}
\begin{document}
\author{Mohamed BOUALI}
\maketitle
\begin{center}{\bf Abstact}
   \end{center}

We describe Generalized Hermitian matrices ensemble sometimes called Chiral ensemble.  We give global asymptotic of the density of eigenvalues or the statistical density. We will calculate a Laplace transform of such a density for finite $n$, which will be expressed through an hypergeometric function. When the dimensional of the hermitian matrix begin large enough, we will prove that the statistical density of eigenvalues converge in the tight topology to some probability measure, which generalize the Wigner semi-circle law.
\section{Introduction}
The eigenvalues distribution of a self-adjoint random $n\times n$ matrix $A$, for $n$ large, was first studied by E. Wigner 1955, and has since been then an active research area in mathematical physics. (See (Mehta 1991) \cite{M} and reference given there).
At my known the model we are aiming it belong to the class of Chiral Models (Akemann 2003 \cite{G}), (Damgaard, Nishigaki 1998, \cite{P}).

In this paper we give an entirely analytical treatment of the key result on asymptotic eigenvalues distribution of generalized gaussian self-adjoint random matrix. All these results generalizes the wigner result. Our ideas is based to derive an explicit formula for the mean value $\mathbb{E}_n\big(Tr[A^2\exp(sA^2)]\big)$ as a function of a complex variable $s$, where $A$ is a self-adjoint random $n\times n$ matrix

We begin by explaining the model, as a model of unitary invariant Hermitian matrices ensemble. In section 2 we give the statistical density (p.d.f) of the eigenvalues and we explain how a scaling of order $\sqrt n$ is necessary for the eigenvalues density. In the same section we recalled some classical results in theory of probability measures, which will be useful for proving the mean theorem. Moreover we will defined the generalized Hermite polynomials which is the crucial tool to derive the Laplace transform of the statistical distribution of eigenvalues. At the end of section 2, we show that, the limit of the statistical density of eigenvalues generalizes the Wigner semi-circle law and it coincides with the Wigner law when some extra parameter become zero.  A such result will be explained on a plot for different values of the extra parameter $c$.

In section 3, by performing some log-Sobolev inequalities, we prove the convergence of the maximal and minimum eigenvalues to the end of the support of the limit measure of the statistical density.

Let $H_n=Herm(n, \Bbb C)$ be the vector space of $n\times n$ Hermitian matrices. For $\mu>-\frac 12$, we denote by $\mathbb{P}_{n,\mu}$
the probability measure on $H_n$ defined by.
$$\int_{H_n}f(x)\mathbb{P}_{n, \mu}(dx)=\frac{1}{C_n}\int_{H_n}f(x)|\det(x)|^{2\mu}e^{-\tr(x^2)}m_n(dx),$$
for a bounded measurable function $f$, where $m_n$ is the Euclidean measure associated to the usual inner product $<x, y>=\tr(xy)$ on $H_n$ and $C_n$ is a normalizing constant.
$$C_n=\prod_{k=0}^{n-1}\alpha_\mu(k)$$
$$\alpha_\mu(k)=\left\{\begin{aligned}&m!\Gamma(m+\mu+\frac12)\quad\mbox {if}\; k=2m,\\
&m!\Gamma(m+\mu+\frac32)\quad\mbox {if}\; k=2m+1.\\
\end{aligned}\right.$$
When $\mu=0$, one recover's the value of the classical Mehta integral $\displaystyle C_n=\pi^{\frac n2}2^{\frac{n(n-1)}{2}}\prod_{k=0}^{n-1}k!$.

We endowed the space $H_n$ with the probability measure $\mathbb{P}_{n, \mu}$. The probability $\mathbb{P}_{n,\mu}$ is invariant for the action of the unitary group $U(n)$ by conjugation
$$x\mapsto uxu^*\qquad (u\in U(n)).$$
\section{Statistical eigenvalues distribution}

Let $f$ be a $U(n)$-invariant function on $H_n$.
$$f(uxu^*)=f(x)\qquad\forall\;u\in U(n),$$
According to the spectral theorem, there exist a symmetric function $F$ in $\Bbb R^n$ such that
$$f(x)=F(\lambda_1,...,\lambda_n).$$
If $f$ is integrable with respect to $\mathbb{P}_{n,\mu}$, then according to the formula of integration of Weyl it follows that
$$\int_{H_n}f(x)\mathbb{P}_{n,\mu}(dx)=\int_{\Bbb R^n}F(\lambda_1,\ldots,\lambda_n)q_n^{\mu}(\lambda_1,...,\lambda_n)d\lambda_1\ldots\lambda_n,$$
where
$$q_n^{\mu}(\lambda_1,\ldots,\lambda_n)=\frac{1}{C_n}e^{-\sum\limits_{k=1}^{n}\lambda_k^2}\prod\limits_{k=1}^{n}\lambda_k^{2\mu}
\Delta(\lambda)^2,$$
and
$$\Delta(\lambda)=\prod_{1\leq i<j\leq n}(\lambda_i-\lambda_j),$$ is the Vandermonde determinant.
We will study the asymptotic of the statistical distribution of the eigenvalues i.e. it is looked as $n$ goes to infinity the asymptotic behavior of the probability $\nu_n$ defined on $\Bbb R$ as follows: if $f$ is a measurable function,
$$\int_{\Bbb R}f(t)\nu_n(dt)=\int_{H_n}\frac1n\tr(f(x))\mathbb{P}_{n,\mu}(dx).$$
Such a measure is absolutely continuous with respect to the Lebesgue measure
$$\nu_n(dt)=h_n^\mu(t)dt.$$
where $$h_n^{\mu}(t)=\int_{\Bbb R^{n-1}}q_n^{\mu}(t,\lambda_2,\ldots,\lambda_n)d\lambda_2\ldots d\lambda_n.$$
We compute the two first moments of the measure $\nu_n$:
$$m_1(\nu_n)=\int_{\Bbb R}t\nu_n(dt)=\frac 1n\int_{H_n}\tr(x)\mathbb{P}_{n,\mu}(dx)=0,$$
the second moment is:
$$m_2(\nu_n)=\frac 1n\int_{H_n}\tr(x^2)\mathbb{P}_{n,\mu}(dx),$$
Since for all $\alpha>0$, $$C_n(\alpha)=\int_{H_n}f(x)|\det(x)|^{2\mu}e^{-\alpha\tr(x^2)}m_n(dx)=\alpha^{-n\mu-n^2-\frac n2}C_n,$$
and $$m_2(\nu_n)=-\frac 1n\frac{d}{d\alpha}\log(C_n(\alpha))|_{\alpha=1}=n+\mu+\frac12.$$
This suggests that $\nu_n$ does't converge, and a scaling of order $\sqrt {n+\mu}$ is necessary. Let defined  $\widetilde \nu_n$ by:
$$\int_{\Bbb R}f(t)\widetilde\nu_n(dt)=\int_{\Bbb R}f(\frac{t}{\sqrt n})\nu_n(dt),$$
We obtain
$$\int_{\Bbb R}f(t)\widetilde\nu_n(dt)=\frac1{Z_n}\int_{\Bbb R^n}\frac1n\sum_{i=1}^nf(x_i)\exp\Big(-n \sum_{i=1}^nQ_{\alpha_n}(x_i)\Big)|\Delta(x)|^\beta dx,$$
where $\displaystyle Q_{\alpha_n}(x)=x^2+\alpha_n\log\frac1{|x|},$ $\alpha_n=\frac{\lambda_n}{n}$, and $Z_n$ is a normalizing constant.

We come to the principal result of the paper
\begin{theorem}--- Let $(\mu_n)_n$ be a sequence of positif real numbers such that $$\lim\limits_{n\rightarrow\infty}\frac{\mu_n}{n}=c\geq 0.$$
Then, after scaling by $\frac{1}{\sqrt n}$, the measure $\nu_n$, converges weakly to the measure
$\nu_c$ supported by $S=[-b,-a]\cup [a, b]$, with density
$$f_c(t)=\frac{1}{\pi |t|}\sqrt{(t^2-a^2)(b^2-t^2)},$$
where $a=\sqrt{1+c-\sqrt {1+2c}}$ and $b=\sqrt{1+c+\sqrt {1+2c}}$. This means that, for a bounded continuous
function $\varphi$ on $\Bbb R$
$$\lim_{n\to\infty}\int_{\Bbb R}\varphi(\frac{u}{\sqrt n})\nu_n(du)=\int_{S}\varphi(u)\nu_c(du).$$
\end{theorem}
Observe that, if $c = 0$, and hence $a = 0, b =\sqrt 2$, one recovers Wigner's
Theorem:
$$\lim_{n\to\infty}\int_{\Bbb R}\varphi(\frac{u}{\sqrt n})\nu_n(du)=\frac1{\pi}\int_{-\sqrt 2}^{\sqrt 2}\varphi(u)\sqrt{2-t^2}du.$$
 We will give a proof by making use the orthogonal polynomials methods.

In forthcoming paper we will study the same model where the Dyson index $2$ in the Vandermonde determinant will be replaced by a strictly positif parameter $\beta >0$.

\subsection{Orthogonal polynomials method}

For $\mu>-\frac12$, the generalized Hermite polynomials $H_n^{\mu}$ are defined by the Rodrigues formula
$$H_n^{\mu}(x)=(-1)^n\frac{n!}{\gamma_{\mu}(n)}e^{x^2}D^n_\mu(e^{-x^2})$$
where $D^n_\mu=D_\mu\circ\cdots\circ D_\mu$ $n$-times, and $D_\mu$ is the one-dimensional Dunkl operator defined by
$$D_{\mu}(\varphi)(x)=\varphi'(x)+\frac{\mu}{x}(\varphi(x)-\varphi(-x)),$$
and
$$\gamma_\mu(k)=\left\{\begin{aligned}&2^{2m}m!\frac{\Gamma(m+\mu+\frac12)}{\Gamma(\mu+\frac12)}\quad\mbox {if}\; k=2m,\\
&2^{2m+1}m!\frac{\Gamma(m+\mu+\frac32)}{\Gamma(\mu+\frac12)}\quad\mbox {if}\; k=2m+1.\\
\end{aligned}\right.$$
See for such definition of generalized Hermite polynomials (Rosenblum 1994 \cite{R}).

The generalized Hermite polynomials are orthogonal with respect to the inner product
$$(p,q)=\int_{\Bbb R}p(x)q(x)|x|^{2\mu}e^{-x^2}dx.$$
In such a way
the square norm of $H_n^{\mu}$ is
$$d_{n,\mu}= \int_{\Bbb R}\big(H_n^{\mu}(x)\big)^2|x|^{2\mu}e^{-x^2}dx=\Gamma(\mu+\frac{1}{2})\frac{2^n(n!)^2}{\gamma_\mu(n)}.$$

We defined the generalized Hermite function as follows
$$\psi_n^{\mu}(x)=\frac{1}{\sqrt {d_{n,\mu}}}|x|^{\mu} e^{-\frac12 x^2}H_n^{\mu}(x).$$
The generalized Hermite polynomials are closely related to the generalized
Laguerre polynomials:
$$H_{2m}^{\mu}(x)=(-1)^m(2m)!\frac{\Gamma(\mu+\frac12)}{\Gamma(m+\mu+\frac12)}L_m^{\mu-\frac12}(x^2),$$
$$H_{2m+1}^{\mu}(x)=(-1)^m(2m+1)!\frac{\Gamma(\mu+\frac12)}{\Gamma(m+\mu+\frac32)}xL_m^{\mu+\frac12}(x^2).$$
Recall that the generalized Laguerre function $\varphi^{\alpha}_n$ is defined, for $\alpha>-1$, $x>0$ by
$$\varphi^{\alpha}_n(x)=\sqrt{\frac{n!}{\Gamma(n+\alpha+1)}}x^{\frac{\alpha}{2}}e^{-\frac x2} L^{\alpha}_n(x).$$
Hence
$$\psi^{\mu}_{2m}(x)=(-1)^m\sqrt{|x|}\varphi^{\mu-\frac12}_m(x^2),$$
$$\psi^{\mu}_{2m+1}(x)=(-1)^m{\rm sign}(x)\sqrt{|x|}\varphi^{\mu+\frac12}_m(x^2),$$
The sequence $\{\psi_n^{\mu}\}_n$ is a Hilbert basis in $L^2(\Bbb R)$.

Before proving the mean result, we give some preliminary results.

\begin{proposition}---
Let $\sigma, \sigma_1,\sigma_2,\cdots $ be a sequence of probability measures on $\Bbb R$ and let ${\cal C}_0(\Bbb R)$, ${\cal C}_b(\Bbb R)$ denote the space of continuous functions on $\Bbb R$, vanishing at $\pm\infty$, respectively the space of bounded continuous functions on $\Bbb R$. Then the following conditions are equivalent
 \begin{enumerate}
\item $\displaystyle\lim_{n\to\infty}\sigma_n(]-\infty,x])=\sigma(]-\infty,x])=F(x)$ for all $x\in\Bbb R$, point of continuity of $F$.
\item For all $\displaystyle f\in{\cal C}_0(\Bbb R), \lim_{n\to\infty}\int_{\Bbb R}f(x)\sigma_n(dx)=\int_{\Bbb R}f(x)\sigma(dx).$
\item For all $\displaystyle f\in{\cal C}_b(\Bbb R), \lim_{n\to\infty}\int_{\Bbb R}f(x)\sigma_n(dx)=\int_{\Bbb R}f(x)\sigma(dx).$
\item For all $\displaystyle x\in\Bbb R, \lim_{n\to\infty}\int_{\Bbb R}e^{itx}\sigma_n(dx)=\int_{\Bbb R}e^{itx}\sigma(dx).$
\end{enumerate}
\end{proposition}
{\bf Proof.} (\cite{fe} Chapter III theorem 1 and theorem 2, chapter VI theorem 2).
\begin{definition}---
Let $\sigma, \sigma_1,\sigma_2,\cdots $ be a sequence of probability measures on $\Bbb R$. If one of the conditions in the previous proposition hold we say that the measure $\sigma_n$ converge weakly to $\sigma$.
\end{definition}
\begin{corollary}---
By the properties $(1)$ of the previous proposition we deduce for all interval $I$ of $\Bbb R$, $\lim_{n\to\infty}\sigma_n(I)=\sigma(I).$
\end{corollary}

 From the fact that the density $h_n^{\mu_n}$ is even function, we will study the asymptotic of $h_n^ {\mu_n}(t)$ as $n$ go to infinity for $t\geq 0$ and we deduce by symmetry the case $t\leq 0$.\\

\begin{proposition}--- The density $h_n^{\mu}$ of the statistical distribution of eigenvalues
is given by
$$h_n^{\mu}( x)=\frac{|x|}{n}\bigg(\sum_{k=0}^{[\frac{n-1}{2}]}\Big(\varphi_k^{\mu-\frac12}(x^2)\Big)^2+\sum_{k=0}^{[\frac {n-2}{2}]}\Big(\varphi_k^{\mu+\frac12}(x^2)\Big)^2\bigg),$$
where $[x]$ is the greatest integer less or equal to $x$.
\end{proposition}
{\bf Proof.}--- Let $K_n^{\mu}$ be the Christoffel-Darboux kernel
$$K_n^{\mu}(x,y)=\sum_{k=0}^{n-1}\psi_k^{\mu}(x)\psi_k^{\mu}(y),$$
The density $h_n^{\mu}$ of the statistical distribution of eigenvalues $\nu_n$ is given by
$$h_n^{\mu}(t)=\frac 1nK_n^{\mu}(t,t).$$
(See for instance (Faraut 2011 \cite{F}) proposition III.3.2 .) Hence the formula for $h_n^{ \mu}$ follows
from the relations above.\\
We will use the following Laplace transform formula.
\begin{lemma}---For ${\rm Re}(s)<1$,
$$\begin{aligned}&\int_0^{+\infty}\sum_{k=0}^{m}\varphi^{\alpha}_k(x)^2x\exp(sx)dx\\
&=(m+1)(m+\alpha+1)\frac{1}{(1-s)^{2m+\alpha+2}}{}_2F_1(-m, -\alpha,2;s^2).
\end{aligned}$$
where ${}_2F_1$ is the Gauss hypergeometric function.
\end{lemma}
(\cite{H}, Lemma 6.3.)

Let $\tau_n$ denote the measure on $\Bbb R_+$ with density ${\sqrt t}h_n^{\mu}(\sqrt t)$. From the
previous Laplace transform formula we get, for ${\rm Re}(s) < n$,
$$\int_0^{+\infty}e^{\frac {s}{n}}\tau_n(dt)=G_n^{\mu, 1}(s)+G_n^{\mu, 2}(s),$$
with
$$\begin{aligned}G_n^{\mu, 1}(s)&=\frac{(m+1)(m+\mu+\frac12)}{n^2}\\&\frac{1}{(1-\frac{s}{n})^{2m+\mu+\frac32}}\;
{}_2F_1(-m-\mu+\frac12, -m,2;\frac{s^2}{n^2}),\end{aligned}$$\\
$$\begin{aligned}G_n^{\mu, 2}(s)&=\frac{(n-m-1)(n-m+\mu-\frac12)}{n^2}\\&\frac{1}{(1-\frac{s}{n})^{2(n-m)+\mu-\frac32}}\;
{}_2F_1(-n+m-\mu+\frac12, 2-n+m,2;\frac{s^2}{n^2}),\end{aligned}$$
where $m=m_n=[\frac{n-1}{2}].$
\begin{lemma}--- Let $\mu_n$ be a sequence of positif real numbers such that
$$\lim_{n\to\infty}\frac{\mu_n}{n}=c.$$
Then $$\lim_{n\to\infty}G^{\mu_n, 1}_n(s)=\lim_{n\to\infty}G^{\mu_n, 1}_n(s)=\frac 12(\frac12+c)e^{(1+c)s}F(\frac{s}{2}\sqrt{1+2c}),$$
where $F$ is the function defined on $\Bbb C$ by
$$F(z)=\sum\limits_{k=0}^{+\infty}\frac{1}{k!(k+1)!}z^{2k}.$$
The convergence is uniform on every compact set in $\Bbb C$.
\end{lemma}
Observe that the function $G^{\mu, 1}_n$ and $G^{\mu, 2}_n$ are defined for ${\rm Re}(s)< n$.
The function F is a Bessel function: $F(z)={}_0F_1(2;z)$.\\
{\bf Proof.}--- We saw
    $$\begin{aligned}&{}_2F_1(-m-\mu+\frac12, -m,2;\frac{s^2}{n^2})\\
    &=\sum_{j=0}^{+\infty}\frac{(m+\mu_n-\frac12)(m+\mu_n-\frac32)\cdots (m+\mu_n-j+\frac12)m(m-1)\cdots (m-j+1)}{j!(j+1)!}\frac{s^{2j}}{n^{2j}},\\
    &=\sum_{j=0}^{+\infty}\frac{a_j(n)}{j!(j+1)!}s^{2j},\\
    &
    \end{aligned}
    $$
with $$a_j(n)=\frac{(m+\mu_n-\frac12)(m+\mu_n-\frac32)\cdots (m+\mu_n-j+\frac12)m(m-1)\cdots (m-j+1)}{n^{2j}}.$$
Since $$\lim_{n\to\infty}\frac{m}{n}=\lim_{n\to\infty}\frac{[\frac{n-1}{2}]}{n}=\frac 12,\;{\rm and}\lim_{n\to\infty}\frac{\mu_n}{n}=c.$$

Then if we put $\gamma=\sup\limits_{n\geq 1}\frac{m+\mu_n}{n}$, since $\frac{m}{n}\leq 1$ we obtain
$$a_j(n)\leq \gamma^j.$$
Let  $R>0$, then for all $|s|\leq R$,
$$\sum_{j=0}^{+\infty}\frac{a_j(n)}{j!(j+1)!}|s|^{2j}\leq e^{\gamma R^2}.$$
By using the Cauchy integral formulas and dominated convergence theorem, we obtain the convergence of the series on every compact subset of $\Bbb C$.
Furthermore, $$\left|(1-\frac{s}{n})^{2m_n+\mu_n+\frac32}\right|\leq\left(1+\frac{|s|}{n}\right)^{2m+\mu_n+\frac32}\leq e^{(3+\gamma)R}.$$

Put all this together, and applying again the dominated convergence theorem for series, we obtain on every compact subset of $\Bbb C$ that
$$\lim_{n\to\infty}G_n^{\mu_n, 1}(s)=\frac 12(\frac12+c)e^{(1+c)s}
\sum\limits_{k=0}^{+\infty}\frac{1}{k!(k+1)!}\Big(\frac{s\sqrt{(1+2c)}}{2}\Big)^{2k}.$$
The same hold for the function $G_n^{\mu_n, 2}$ if we replace the sequence $a_j(n)$ by
$$b_j(n)=\frac{(n-m+\mu_n-\frac12)\cdots (n-m+\mu_n-j+\frac12)(n-m-1)(n-m-2)\cdots (n-m-j-1)}{n^{2j}}.$$
Since for all $k\in\Bbb N$,
$$n-m-k\leq n\;\;\mbox{and}\;\; n-m-k+\mu_n-\frac12\leq m+\mu_n+\frac52, $$
then $$b_j(n)\leq(\gamma+3)^j,$$
and we deduce the same limits on every compact subset of $\Bbb C$ as for the function $G_n^{\mu_n, 1}$.

Therefore
\begin{equation}\label{1}\begin{aligned}\lim\limits_{n\to +\infty}\int_0^{+\infty}e^{s\frac{t}{n}}\tau_n(dt)&=\lim_{n\to\infty}G_n^{\mu_n, 1}(s)+\lim_{n\to\infty}G_n^{\mu_n, 2}(s)\\&=(\frac12+c)e^{(1+c)s}
{}_0F_1(2;\frac{s}{2}\sqrt{1+2c)}).\end{aligned}\end{equation}

\begin{lemma}--- For $\alpha<\beta$, $s\in\Bbb C$,
$$\frac1{\pi}\int_{\alpha}^{\beta}e^{st}\sqrt{(t-\alpha)(\beta-t)}dt=\frac12\left(\frac{\beta-\alpha}{2}\right)^2e^{\frac{\beta+\alpha}{2}s}{}_0F_1(2;\frac{\beta-\alpha}{4}s).$$

\end{lemma}
{\bf Proof.}--- By performing the change of variables $u=\frac{2}{\beta-\alpha}(t-\frac{\beta+\alpha}{2})$, it follows
$$\frac1{\pi}\int_{\alpha}^{\beta}e^{st}\sqrt{(t-\alpha)(\beta-t)}dt=(\frac{\beta-\alpha}{2})^2e^{\frac{\alpha+\beta}{2}s}\frac1{\pi}\int_{-1}^1e^{\frac{\beta-\alpha}{2}s}\sqrt{1-u^2}du,$$
The last integral is essentially a Bessel function and the result follow.\\

{\bf Proof of theorem 2.1.}--- Let denote by $g_c$ the density on $[a^2,b^2]$ defined by
$$g_c(t)=\frac{1}{\pi }\sqrt{(t-a^2)(b^2-t)},$$
where $a^2=1+c-\sqrt {1+2c}$ and $b^2=1+c+\sqrt {1+2c}$.

From the previous lemma we have for all $s\in\Bbb C$, \begin{equation}\label{2}\int_0^{+\infty}e^{st}g_c(dt)=(\frac12+c)e^{(1+c)s}{}_0F_1(2;\frac{s}{2}\sqrt{1+2c})
.\end{equation}
Since the function ${}_0F_1$ is continuous at $0$, applying L\'evy-Cram\'er theorem we deduce from equations (\ref{1}) and (\ref{2}) that, the measure $\tau_n$ scaled by $\frac1n$ converge weakly to the measure $\tau$ with density $g_c$ with respect to the Lebesgue measure .

Since $$\tau_n(dt)=\sqrt th_n^{\mu_n}(\sqrt t)dt,$$
hence, for every compactly supported function $\varphi\in{\cal C}_c(\Bbb R^+)$ and by performing the change of variable $t=\sqrt u$ we get,
$$\begin{aligned}\int_0^{+\infty}t^2\varphi(t)h_n^{\mu_n}(t)dt&=\frac12\int_0^{+\infty}\varphi(\sqrt u)\sqrt uh_n^{\mu_n}(\sqrt u)du\\
&=\frac12\int_0^{+\infty}\varphi(\sqrt u)\tau_n(du).
\end{aligned}$$

Moreover $$\lim_{n\to +\infty}\frac12\int_0^{+\infty}\varphi(\sqrt{\frac{ u}{n}})\tau_n(du)=\frac12\int_{a^2}^{b^2}\varphi(\sqrt u)g_c(u)du=\int_{ a}^{ b}\varphi( v)g_c(v^2)vdv.$$
it follows that,
$$\lim_{n\to +\infty}\int_0^{+\infty}(\frac{t}{\sqrt n})^2\varphi(\frac{t}{\sqrt n})h_n^{\mu_n}(t)dt=\int_{a}^{ b}\varphi( t)g_c(t^2)tdt.$$

Therefore, there is some constant $A\geq 0$, such that for every $\varphi\in{\cal C}_c([0,+\infty[)$,
$$\lim_{n\to +\infty}\int_0^{+\infty}\varphi(\frac t{\sqrt n})h_n^{\mu_n}(t)dt=A\varphi(0)+\int_{a}^{ b}\varphi(t)g_c(t^2)\frac{dt}{t}.$$
Since the function $h_n^{\mu_n}$ is even, It readily follows that, for all $\psi\in{\cal C}_c(\Bbb R)$,
$$\lim_{n\to +\infty}\int_{\Bbb R}\psi(\frac t{\sqrt n})h_n^{\mu_n}(t)dt=A\psi(0)+\int_{[-b,-a]\cup[a,b]}\psi(t)g_c(t^2)\frac{dt}{|t|}.$$
To compute the constant $A$, observe that $$\int_{[-b,-a]\cup[a,b]}g_c(t^2)\frac{dt}{|t|}=1.$$
Hence the limit measure is a probability and $A$=0.

Which means that after scaling by $\frac{1}{\sqrt n}$ the density $h_n^{\mu_n}$ converge as $n$ goes to infinity on every compact subset of $\Bbb R$ to the function $f_c$ where $$f_c(t)=\frac{1}{\pi |t|}\sqrt{(t^2-a^2)(b^2-t^2)},$$
 $a=\sqrt{1+c-\sqrt {1+2c}}$ and $b=\sqrt{1+c+\sqrt {1+2c}}$.\\
This complete the proof.

\begin{figure}[h]
\centering{\scalebox{0.75}{\includegraphics[width=15cm, height=6cm]{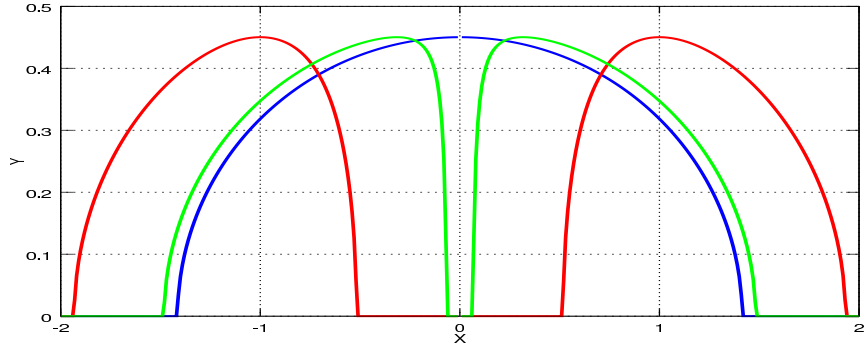}}}
\caption{\bf Plot of the density $f_c$.}
\end{figure}
\begin{flushleft}\textcolor{green}{\rule{1cm}{1.5pt}} for $c=0.1$.\\
\textcolor{red}{\rule{1cm}{1.5pt}} for $c=1$.\\
\textcolor{blue}{\rule{1cm}{1.5pt}}\,for $c=0$, the semi-circle law.\\
\end{flushleft}

\section{Asymptotic of the largest eigenvalue}
\begin{proposition}---
Let $X_n$ be sequence of hermitian random matrices defined on the same probability space and such that the distribution of $X_n$ is the probability $\mathbb{P}_{n, \mu_n}$ for all $n\in\Bbb N$. Let $\mu_n$ be a positive real sequence such that $\displaystyle\lim_{n\to\infty}\frac{\mu_n}{n}=c$ for some constant $c$, then
\begin{enumerate}
\item$\lim_{n\to +\infty}\frac{\lambda_{\max}(X_n)}{\sqrt n}= b\quad{\rm almost\; surely}.$
\item$\lim_{n\to +\infty}\frac{\lambda_{\min}(X_n)}{\sqrt n}= -b\quad{\rm almost\; surely}.$
\end{enumerate}

\end{proposition}
Before proving the proposition we will give sum preliminary results in concentration inequalities.
\subsection{ Concentration inequalities}
\begin{definition}---A probability measure $\sigma$ on $\Bbb R^n$ is said to be satisfy a logarithmic Sobolev inequality if there exists some nonnegative constant $c_0$ independent of $n$ such that,
 \begin{equation}\int_{\Bbb R^n}f^2\log(f^2)d\sigma\leq 2c_0\int_{\Bbb R^n} ||\nabla f||_2^2d\sigma,\end{equation}
 for every smooth function $f:\Bbb R^n\longrightarrow\Bbb R$, such that $\displaystyle\int_\Bbb Rf^2d\mu=1$.

 Where $\nabla$ is the Euclidean gradient and $||\cdot||_2$ the Euclidean norm on $\Bbb R^n$.
 \end{definition}

 The next proposition gives an important example related to Boltzmann density.
 \begin{proposition}{\rm (Ledoux 1999, \cite{Mi})}

Let $\nu$ be a Boltzmann probability on $\Bbb R^n$ of type $d\mu=\frac{1}{Z}\exp(-V(x))dx,$ such that $V$ is a twice continuous differentiable function. If there is some constant $c>0$ such that $\displaystyle V''(x)-cI_n$ is positive symmetric matrix for all $x\in\Bbb R^n$, then $\mu$ satisfies $(3.3)$ with constant $\displaystyle c_0=\frac{1}{c}$

 \end{proposition}

 \begin{corollary} {\rm ( Ledoux 1999, \cite{Mi})}

 Let $\nu$ be a Boltzmann measure on $\Bbb R^n$ satisfying equation $(3.3)$ and $F$ be a real or complex valued Lipschitz function on $\Bbb R^n$ such that $||F||_{Lip}\leq\alpha$,
 then

 $$\nu\Big\{x\in\Bbb R^n,\, |F(x)-E_{\nu}(F)|>r\Big\}\leq 2\exp(-r^2/{4c_0\alpha^2}).$$
 Where $E_\nu$ is the expectation of $F$ with respect to the measure $\nu$ and $\displaystyle||F||_{Lip}=\sup_{x\in\Bbb R^n}\frac{|F(x)|}{||x||_2}$ is the Lipschitz constant.
 \end{corollary}
\begin{proposition}---
Let $\gamma_n$ be the probability measure defined on the open set $$\Sigma=\{(x_1,\cdots,x_n)\in\Bbb R^n\mid x_i\neq x_j,\;x_i\neq 0,\,\forall\,1\leq i<j\leq n\},$$  with density $\displaystyle\frac{1}{Z_n}e^{-V_n(x)},$
  where $$\displaystyle V_n(x)=\sum_{i=1}^n(x_i^2-2\mu_n\log|x_i|)-2\sum_{1\leq i<j\leq n}\log|x_i-x_j|,$$  $\mu_n$ is a sequence of positif real numbers and $Z_n$ is a normalizing constant. Then $\gamma_n$  satisfies $(3.3)$ with constant $\displaystyle c_0=\frac12$.\\

\end{proposition}
\begin{corollary}--- Let $f$ be a Lipschitz function on $\Bbb R$ with $||f||_{Lip}\leq\alpha$. For $\delta>0$, defined on $\Bbb R^n$ the function
$$F(\frac{x}{\delta})=\frac1n\sum_{i=1}^nf(\frac{x_i}{\delta})\quad \forall\;x=(x_1,\cdots, x_n)\in\Bbb R^n.$$ Then
\begin{enumerate}
\item The function $F$ is lipschitz on $\Bbb R^n$ with $\displaystyle ||F||_{Lip}\leq\frac{\alpha}{\delta\sqrt n}.$
\item For all $\e>0$, and for each $\delta>0$, $$\gamma_n\left\{x\in\Bbb R^n,\; \Big|F(\frac{x}{\delta})-E_{\gamma_n}(F)\Big|>\e\right\}\leq 2 \exp(-\frac{n\delta^2\e^2}{2{\alpha^2}}),$$
    where $E_{\gamma_n}$ is the expectation with respect to the probability $\gamma_n$.
 \end{enumerate}
\end{corollary}
{\bf Proof.}---\\
1) We know that $||f||_{Lip}\leq\alpha$, which means that for all $s,t\in\Bbb R$,
$$|f(s)-f(t)|\leq\alpha |s-t|.$$
Hence by the Cauchy-Schwartz inequality it follows that, for all $\delta>0$
$$\Big|F(\frac x{\delta})-F(\frac y\delta)\Big|\leq \frac1n \sum_{i=1}^n|f(\frac {x_i}{\delta})-f(\frac {y_i}{\delta})|\leq\frac{\alpha}{\delta\sqrt n}||x-y||_2,$$
where $||\cdot||_2$ is the Euclidean norm on $\Bbb R^n$. Hence
$$||F||_{Lip}\leq\frac{\alpha}{\delta\sqrt n}.$$
2) The function $V_n$ is twice differentiable on the open set \\$\Big\{(x_1,\cdots,x_n)\in\Bbb R^n, x_i\neq x_j,\;x_i\neq 0\Big\}$. Furthermore
for all $k=1,\cdots ,n$
$$\frac{\partial V_n}{\partial x_k}=2x_k-\frac{2\mu_n}{x_k}-2\sum_{i=1, i\neq k}^n\frac{1}{x_k-x_i}.$$
Take the second derivative, it follows
$$\frac{\partial^2 V_n}{\partial^2 x_k}=2+\frac{2\mu_n}{x_k^2}+2\sum_{i=1, i\neq k}^n\frac{1}{(x_k-x_i)^2},$$
and for $\ell\neq k$,
$$\frac{\partial^2 V_n}{\partial x_\ell\partial x_k}=\frac{-2}{(x_k-x_\ell)^2}.$$
 The goal is to  prove that $\displaystyle V''(x)-cI_n$ is  positive matrix uniformly on $x$ in the sense of the usual Euclidean structure .\\
Let $y=(y_1,\cdots,y_n)$ be an element in $\Bbb R^n$.
$$<V_n''(x)y,y>=2||y||^2+2\mu_n\sum_{k=1}^n\frac{y_k^2}{x_k^2}+2\sum_{i\neq k}\frac{y^2_k}{(x_k-x_i)^2}-2\sum_{i\neq k}\frac{y_ky_i}{(x_k-x_i)^2}.$$
Now applied  two time the Cauchy-Schwartz inequality, this gives
$$\begin{aligned}&\sum_{i\neq k}\frac{y_ky_i}{(x_k-x_i)^2}\\
&\leq\sum_{i=1, i\neq k}^{n}\Big[\Big(\sum_{k=1}^{n}\frac{y_k^2}{(x_k-x_i)^2}\Big)^{\frac{1}{2}}\Big(\sum_{k=1}^{n}\frac{y_i^2}{(x_k-x_i)^2}\Big)^{\frac{1}{2}}\Big]\\
&\leq\Big(\sum_{i=1, i\neq k}^{n}\sum_{k=1}^{n}\frac{y_k^2}{(x_k-x_i)^2}\Big)^{\frac{1}{2}}\Big(\sum_{i=1, i\neq k}^{n}\sum_{k=1}^{n}\frac{y_i^2}{(x_k-x_i)^2}\Big)^{\frac{1}{2}}\\
&\leq \sum_{i\neq k}\frac{y^2_k}{(x_k-x_i)^2}.
\end{aligned}$$
Hence $$<V''(x)y,y>\;\geq 2||y||^2\quad\forall\;x\in\Sigma,$$
which means that the symmetric matrix $V''(x)-2I_n$ is positive.\\

We come to our application
\begin{proposition}---
Let $X_n$ be a sequence of hermitian random matrices, defined on the same probability space $(\Omega,{\cal B},{\cal P})$, such that $X_n$ follows as the probability law $\displaystyle\mathbb{P}_{n, \mu_n}(dX)=\frac{1}{C_n}|\det(X)|^{2\mu_n}e^{-\tr(X^2)}m_n(dX),$ where $C_n$ is a normalizing constant, $\mu_n$ is a positive sequence and $m_n$ is the Lebesgue measure on the space $H_n$ of hermitian matrices.

For each $n\in\Bbb N$, and for every $\omega\in\Omega$, let $\rho_{n,\omega}$ denoted the statistical distribution of the ordering eigenvalues $\lambda_1\big(\frac{X_n(\omega)}{\sqrt n}\big)\leq\lambda_2\big(\frac{X_n(\omega)}{\sqrt n}\big)\leq\cdots\leq\lambda_n\big(\frac{X_n\big(\omega)}{\sqrt n}\big)$ of $X_n(\omega)$, which means that
$$\rho_{n,\omega}=\frac1n\sum_{i=1}^n\delta_{\lambda_i\big(\frac{X_n(\omega)}{\sqrt n}\big)}.$$ Then for almost all $\omega\in\Omega$, under the condition $\displaystyle\lim\limits_{n\to\infty}\frac{\mu_n}{n}=c$, the measure $\rho_{n,\omega}$ converge weakly to the distribution $\nu_c$, with density
$$f_c(t)=\frac{1}{\pi |t|}\sqrt{(t^2-a^2)(b^2-t^2)}\chi_S(t),$$
where $\displaystyle S=[-b,-a]\cup [a,b]$, $\displaystyle a=\sqrt{1+c-\sqrt {1+2c}}$, $b=\sqrt{1+c+\sqrt {1+2c}}$ and $\chi_S$ is the characteristic function of the set $S$.\\
Furthermore, for almost all $\omega\in\Omega$, and all interval $I$ in $\Bbb R$
$$\lim_{n\to\infty} \left(\frac1n{\rm card}\big({\rm sp}[X_n(\omega)]\big)\cap I\right)=\nu_c(I),$$
where ${\rm sp}[X_n(\omega)]$ mean the spectrum of the random matrix $X_n(\omega)$.
\end{proposition}
{\bf Proof.}--- In the first hand, for all $f\in{\cal C}_0(\Bbb R)$,.\\
$$\int_{\Bbb R}f(x)\rho_{n,\omega}(dx)=\frac1n\sum_{i=1}^nf\left(\lambda_i(\frac{X_n(\omega)}{\sqrt n})\right)=\tr_n\left[f\left(\frac{X_n(\omega)}{\sqrt n}\right)\right],$$
where ${\cal C}_0(\Bbb R)$ is the space of continuous functions which go to $0$ at infinity.\\
So it suffices to show that,

$$\lim_{n\to\infty}\tr_n\left[f\left(\frac{X_n(\omega)}{\sqrt n}\right)\right]=\int_{\Bbb R}f(x)\nu_c(dx),\qquad {\rm for\; almost\; all}\; \omega\in\Omega\; {\rm and\; all\;}f\in{\cal C}_0(\Bbb R).$$
By separability of the Banach space ${\cal C}_0(\Bbb R)$, one can just prove the previous equality for each $f$ in some dense subspace of ${\cal C}_0(\Bbb R)$, the space ${\cal C}^1_c(\Bbb R)$ of continuous differentiable functions on $\Bbb R$ with compact support. Let $f$ be a function in ${\cal C}^1_c(\Bbb R)$ and put
$$F(X)=\tr_n(f(X)),\;\forall\,X\in H_n,$$
The function  $f$ is Lipschitz with constant $\alpha=||f||_{Lip}=\sup_{x\in\Bbb R}|f'(x)|$, hence by corollaty 3.6, the function $F$ is Lipschitz too with $||F||_{Lip}\leq\frac{\alpha}{\sqrt n}$ .\\
 The function $F$ and the probability $\displaystyle\mathbb{P}_{n, \mu_n}$ are invariants by conjugation by the action of the unitary group. Hence if we denote by $\theta_n$, the invariant part of the measure $\displaystyle\mathbb{P}_{n, \mu_n}$, then for all continuous invariant function $\varphi$ on $H_n$,
 $$\int_{H_n}\varphi(X) \mathbb{P}_{n, \mu_n}(dX)= \int_{\Bbb R^n}\phi(\lambda_1,...,\lambda_n) \theta_n(d\lambda_1,...,d\lambda_n),$$
 where $\varphi(X)=\phi(\lambda_1,...,\lambda_n)$, and $\theta_n$ is the probability on $\Bbb R^n$ with density
 $$\frac1{C_n}e^{-\sum\limits_{i=1}^n\lambda_i^2}\prod_{i=1}^n|\lambda_i|^{2\mu_n}\prod_{1\leq i<j\leq n}(\lambda_i-\lambda_j)^2.$$
 Let $\mathbb{E}_n$, respectively $E_{\theta_n}$ be the expectation with respect the probability  $\displaystyle\mathbb{P}_{n, \mu_n}$, respectively $\theta_n$.

 Observe that, for all $\e >0$, the set $$\{X_n(\omega),\; |F(X_n)-\displaystyle\mathbb{E}_n(F)|>\e\}$$ is invariant by conjugation.

 Therefore by the spectral theorem
 $${\cal P}\big\{\omega\in\Omega,\;|F\big(\frac{X_n(\omega)}{\sqrt n}\big)-\displaystyle\mathbb{E}_n(F)|>\e\big\}={\cal P}\big\{\omega\in\Omega,\;|F\big(\frac{x_n(\omega)}{\sqrt n}\big)-\displaystyle\mathbb{E}_{\theta_n}(F)|>\e\big\},$$
where $x_n(\omega)=\Big(\lambda_1(X_n(\omega)),\lambda_2(X_n(\omega)),\cdots,\lambda_n(X_n(\omega))\Big)\in\Bbb R^n$.

Furthermore $$\Bbb R^n\setminus\Sigma =\bigcup_{ i,j =1}^nH_{i,j},$$
where $H_{i,j}$ is the vector subspace $\displaystyle H_{i,j}=\{x\in\Bbb R^n\mid x_i=x_j, i\neq j\}$ for $i\neq j$, and $H_{i,i}=\{x\in\Bbb R^n\mid x_i=0, 1\leq i\leq n\}$ which is an hyperplane. Since for all $i, j$, $\theta_n(H_{i,j})=0$, hence
$\theta_n(\Bbb R^n\setminus\Sigma)=0.$ Which means that the two measures $\theta_n$ and $\gamma_n$ are equal on $\Sigma$ where $\gamma_n$ is the measure
of the proposition\,3.5.
Hence
 $${\cal P}\big\{\omega\in\Omega,\;|F\big(\frac{X_n(\omega)}{\sqrt n}\big)-\displaystyle\mathbb{E}_n(F)|>\e\big\}={\cal P}\big\{\omega\in\Omega,\;|F\big(\frac{x_n(\omega)}{\sqrt n}\big)-E_{\gamma_n}(F)|>\e\big\},$$
 here $E_{\gamma_n}$ means the expectation with respect to the probability $\gamma_n$ and $x_n(\omega)\in\Sigma$.

From corollary\,3.6,  we deduce that
$${\cal P}\Big\{\omega\in\Omega,\;|F\big(\frac{x_n(\omega)}{\sqrt n}\big)-E_{\gamma_n}(F)|>\e\Big\}\leq 2\exp(-\frac{n^2\e^2}{2\alpha^2}).$$
By using the Borel Cantelli lemma, it follows that, for almost all $\omega\in\Omega$
 $$|F\big(\frac{x_n(\omega)}{\sqrt n}\big)-\displaystyle\mathbb{E}_n(F)|\leq\e,\qquad\;\forall\,\e>0.$$
 Since the space ${\cal C}^1_c(\Bbb R)$ is dense in ${\cal C}_0(\Bbb R)$, hence
 $$\lim_{n\to\infty}F\big(\frac{x_n(\omega)}{\sqrt n}\big)=\lim_{n\to\infty}\mathbb{E}_n(F)=\int_{\Bbb R}f(x)\nu_c(dx),\qquad\forall\;f\in{\cal C}_0(\Bbb R).$$
 Which mean that
 $$\lim_{n\to\infty}\int_{\Bbb R}f(x)\rho_{n,\omega}(dx)=\int_{\Bbb R}f(x)\nu_c(dx),\qquad{\rm for\,all}\;f\in{\cal C}_0(\Bbb R), {\rm and\, all}\, \omega\in\Omega.$$
 The second step of the proposition follow from a classical results of convergence of measures theory.

 {\bf Proof of proposition 3.1.}--- Now we will show the convergence of the largest eigenvalue.\\

    {\bf First step:}
 We first show that $$\limsup_{n\to +\infty}\frac{\lambda_{\max}(X_n)}{\sqrt n}\leq b.$$
Let $X_n$ be a random matrix follows as the distributed $\mathbb{P}_{n,\mu_n}$, and put
$$b_n^2=\frac{n+\mu_n+2+\sqrt{n(n+2\mu_n+4)}}{n}.$$
Observe that the sequence $b_n^2$ converge to $b^2$.

Moreover, for all $t\geq 0$, let $\psi_{n,X_n}$ by defined as follows
$$\psi_{X_n}(t)=\mathbb{E}_n(X_n^2\tr(\exp(tX_n^2))).$$
where $\mathbb{E}_n$ is the expectation with respect to the probability $\mathbb{P}_{n,\mu_n}$.

In first hand
$$\psi_{X_n}(t)=n\int_{\Bbb R}x^2e^{tx^2}h_n^{\mu_n}(x)dx
=n\int_{0}^{+\infty}e^{tx}\tau_{n}(dx)
=n\big(G_n^{\mu_n, 1}(nt)
+G_n^{\mu_n, 2}(nt)\big),
$$

furthermore, for all $t\in[0,\frac 12]$,
$$\frac{1}{1-t}\leq e^{t+t^2},\;{\rm and }\,{}_2F_1(a, b,2;t^2)\leq e^{2t\sqrt{|ab|}}.$$

Thus according to the definition of the functions $G_n^{\mu_n, 1}$,
$G_n^{\mu_n, 2}$, one gets
$$\begin{aligned}&\psi_{n, X_n}(t)\\&\leq M_n\left(e^{(2m+\mu_n+\frac32)(t+t^2)+2t\sqrt{(m+\mu_n-\frac12)(m-2)}}+e^{(2(n-m)+\mu_n-\frac32)(t+t^2)+2t\sqrt{(n-m+\mu_n-\frac32)(n-m-2)}}\right),\end{aligned}$$
where $M_n=\max(a_n,b_n),$ $\displaystyle a_n=\frac{(m+1)(m+\mu_n+\frac12)}{n}$ and $\displaystyle b_n=\frac{(n-m-1)(n-m+\mu_n-\frac12)}{n}$.\\\\
Furthermore by simple computation we have,
$$\psi_{n, X_n}(t)\leq\delta_n e^{nb_n^2 t+(n+\mu_n+4)t^2}.$$
 Let $Y_n$ be the random matrix defined by $Y_n=\frac {1}{\sqrt n}X_n$, then
 $$\psi_{n, Y_n}(t)=\frac 1n\psi_{n, X_n}(\frac{t}{n}),$$
 and for $t\in[0,\frac n2]$,
 $$\psi_{n, Y_n}(t)\leq\frac{\delta_n}{n}e^{b_n^2t+\frac 1n(1+\frac{\mu_n}{n}+\frac 4n)t^2}.$$
 Furthermore $Y_n$ and $Y_n^2$ are hermitian matrices, by spectral theorem and straightforward computation we obtain
$$\Big(\lambda_{\max}(Y_n)\Big)^2\leq \lambda_{\max}(Y_n^2).$$
Thus for all $t\geq 0$, by positivity of the exponential function we have
$$\tr\Big(\exp(tY_n^2)\Big)\geq\exp\Big(t\lambda_{\max}(Y_n^2)\Big)\geq\exp\Big(t\big(\lambda_{\max}(Y_n)\big)^2\Big),$$
then,
$$\begin{aligned}\mathbb{P}_{n,\mu_n}\Big(\big(\lambda_{\max}(Y_n)\big)^2\geq b_n^2+\e\Big)&=\mathbb{P}_{n,\mu_n}\Big(t\Big\{\big[\big(\lambda_{\max}(Y_n)\big)^2- \big(b_n^2+\e\big)\big]\Big\}\geq 0\Big)\\
&=\mathbb{P}_{n,\mu_n}\bigg(\exp\Big(t\Big[\big(\lambda_{\max}(Y_n)\big)^2- \big(b_n^2+\e\big)\Big]\Big)\geq 1\bigg)\\
&\leq\mathbb{P}_{n,\mu_n}\Big(\exp\big(t\lambda_{\max}(Y_n^2)- t(b_n^2+\e)\big)\geq 1\Big)\\
&\leq \mathbb{E}_n\Big\{\exp\big(- t(b_n^2+\e)\big)Y_n^2\exp\big(tY_n^2\Big)\Big\}\\
&=e^{-t(b_n^2+\e)}\mathbb{E}_n\left(Y_n^2\exp(tY_n^2)\right)\\
&=e^{-t(b_n^2+\e)}\psi_{n,Y_n}(t)\\
&\leq\frac{\delta_n}{n} e^{-t\e+(\frac 1n+\frac{\mu_n}{n}+\frac 4n)\frac{t^2}{n}}\\
&\leq \frac{\delta_n}{n}e^{-t_0\e+(\frac 1n+\frac{\mu_n}{n}+\frac 4n)\frac{t_0^2}{n}},
\end{aligned}$$
where $$t_0=\frac{n\e}{2(\frac 1n+\frac{\mu_n}{n}+\frac 4n)}.$$
 We have used the fact that the function $t\mapsto -t\e+(\frac 1n+\frac{\mu_n}{n}+\frac 4n)\frac{t^2}{n}$ attaint it minimum at the point $t_0=\frac{n\e}{2(\frac 1n+\frac{\mu_n}{n}+\frac 4n)}$.

It follows that
$$\mathbb{P}_{n,\mu_n}\Big(\big(\lambda_{\max}(Y_n)\big)^2\geq b_n^2+\e\Big)\leq\frac{\delta_n}{n}\exp\Big({\frac{-n\e}{4(\frac 1n+\frac{\mu_n}{n}+\frac 4n)}}\Big).$$
Observe that the sequence $\displaystyle\frac{\delta_n}{n}=\sup\Big(\frac{(m+1)(m_n+\mu_n+\frac32)}{n};\frac{(n-m-1)(n-m+\mu_n-\frac12)}{n}\Big)$ is bounded.\\
Furthermore as $n\to +\infty$, $$\exp\left({\frac{-n\e}{4(\frac 1n+\frac{\mu_n}{n}+\frac 2n)}}\right)\sim e^{-\frac{n\e}{4c}}.$$
This implies that
$$\sum_{n=1}^{\infty}\mathbb{P}_{n,\mu_n}\Big(\big(\lambda_{\max}(Y_n)\big)^2\geq b_n^2+\e\Big)<\infty.$$
Hence the Borel-Cantelli lemma yield, that on a set with probability one, we have
$$\Big(\lambda_{\max}(Y_n)\Big)^2\leq b_n^2+\e,\quad \mbox{eventually},$$
and consequently
$$\limsup_{n\to +\infty}\Big(\lambda_{\max}(Y_n)\Big)^2\leq\limsup_{n\to +\infty} b_n^2+\e\leq b^2+\e.$$
 Since $b\geq 0$, hence $b^2+\e\leq (b+\sqrt\e)^2$, thus for all $\e> 0$,
$$\limsup_{n\to +\infty}\lambda_{\max}(Y_n)\leq b+\e,$$
and
$$\limsup_{n\to +\infty}\lambda_{\max}(Y_n)\leq b.$$
Which means that, if $X_n$ is a random matrix distributed according to the probability $\mathbb{P}_{n,\mu_n}$ then
$$\limsup_{n\to +\infty}\frac{\lambda_{\max}(X_n)}{\sqrt n}=\limsup_{n\to +\infty}\lambda_{\max}(\frac{X_n}{\sqrt n})=\limsup_{n\to +\infty}\lambda_{\max}(Y_n)\leq b.$$

{\bf Second step.}
  In the other hand given any $\e>0$, from the previous proposition, for almost all $\omega\in\Omega$,
  $$\lim_{n\to\infty} \left(\frac1n{\rm card}({\rm sp}[X_n(\omega)])\cap [b-\e,+\infty[\right)=+\infty.$$
  Which mean that
   $\displaystyle\liminf_{n\to +\infty}\frac{\lambda_{\max}(X_n)}{\sqrt n}\geq b-\e\quad{\rm for \,almost\, all}\; \omega\in\Omega,\,{\rm and\, for\, every\,\e>0},$
   then the conclusion hold.
    $$\displaystyle\lim_{n\to +\infty}\frac{\lambda_{\max}(X_n)}{\sqrt n}=b\quad{\rm for \,almost\, all}\; \omega\in\Omega.$$
    The second step in the proposition follows by considering $-X_n$ instead of $X_n$.
\begin{center}{\bf Acknowledgments.}
\end{center}
My sincere thanks go to Jacques Faraut for his comments on this manuscript and his important remarks. Thanks also to the referee for here useful remarks.

Department of Mathematics, College of Applied Sciences,
Umm Al-Qura University, P.O Box  (715), Makkah,
Saudi Arabia.\\
 E-mail address: bouali@math.jussieu.fr, bouali25@laposte.net\\

\end{document}